\documentclass[12pt]{article}
\usepackage[T2A]{fontenc}
\usepackage[utf8]{inputenc}
\usepackage{graphicx,xcolor}
\usepackage{amssymb,amsmath,amsfonts,amsthm,amscd,latexsym,verbatim,graphics,epsfig,indentfirst,color}
\usepackage{geometry}
\def\id{\mathrm{id}}
\def\Aut{\mathrm{Aut}}
\def\Spec{\mathrm{Spec}}
\def\Imm{\mathrm{Im}\,}

\begin{document}

\hfill{16W99, 17C20 (MSC2020)}

\begin{center}
{\Large
Rota---Baxter operators on the simple Jordan algebra of matrices of order two}

V. Gubarev, A. Panasenko
\end{center}

\begin{abstract}
We describe all Rota---Baxter operators of any weight on the space of matrices from $M_2(F)$ considered under the product $a\circ b = (ab + ba)/2$ and usually denoted as $M_2(F)^{(+)}$.
This algebra is known to be a simple Jordan one.

We introduce symmetrized Rota---Baxter operators of weight~$\lambda$
and show that every Rota---Baxter operator of weight~0 on $M_2(F)^{(+)}$
either is a Rota---Baxter operator of weight~0 on $M_2(F)$ or is a~symmetrized Rota---Baxter operator of weight~0 on the same $M_2(F)$.

We also prove that every Rota---Baxter operator of nonzero weight~$\lambda$ on $M_2(F)^{(+)}$
is either a Rota---Baxter operator of weight~$\lambda$ on $M_2(F)$ or 
is, up to the action of $\phi\colon R\to -R-\lambda\id$, a~symmetrized Rota---Baxter operator of weight~$\lambda$ on $M_2(F)$. 

\end{abstract}

\section{Introduction}

Let $A$ be an algebra defined over the ground field~$F$ and fix $\lambda\in F$.
A~linear operator $R\colon A\rightarrow A$ is called a Rota---Baxter operator
(RB-operator, for short) of weight~$\lambda$ if the identity
\begin{equation}\label{RB}
R(x)R(y) = R( R(x)y + xR(y) + \lambda xy )
\end{equation}
holds for all $x,y\in A$.
An algebra~$A$ equipped with a~Rota---Baxter operator is called
a~Rota---Baxter algebra (RB-algebra, for short).

Rota---Baxter operators were defined by G. Baxter in 1960~\cite{Baxter}.
When $\lambda = 0$, the relation~\eqref{RB} may be considered
as an algebraic analogue of the integration by parts formula.
In 1980s, such operators were rediscovered in the context of different versions of Yang---Baxter equation~\cite{Semenov83}.

In~2000, M. Aguiar showed~\cite{Aguiar} that solutions
of the associative Yang---Baxter equation and Rota---Baxter operators of weight zero on any associative algebra are related.
In~2000, a~deep connection between RB-algebras and so-called prealgebras was also found~\cite{Aguiar-2}.
In~2012, Li Guo wrote a monograph devoted to Rota---Baxter algebras~\cite{GuoMonograph}.

One of the directions in the study of RB-operators is a partial or complete classification of them on a given algebra. Since the relation~\eqref{RB} is quadratic, the solution of the problem may not be easy. A special interest lies in the description of RB-operators on simple algebras.
Rota---Baxter operators were classified of any weight on $M_2(F)$ in~\cite{Aguiar,BGP,Mat2} (see also~\cite{GonGub,Panasenko}), of nonzero weight on~$M_3(F)$ in a series of works~\cite{GonGub,Gub2021,Gub2024}, skew-symmetric ones of weight zero on $M_3(F)$ in~\cite{Sokolov},
of any weight on the Cayley---Dickson algebra in~\cite{Panasenko,Panasenko2}, for more details see~\cite{Unital}.

Rota---Baxter operators on Jordan algebras have been studied in a series of works~\cite{post-Jordan,JordRep,preJordan,23dimJord,2dimpreJord}.
However, Rota---Baxter operators on simple Jordan algebras were considered only on spin factors~\cite{Gubarev2017,BGP,Spectrum}.
Recall that any associative algebra $(A,\cdot)$ under the new product
$a\circ b = (a\cdot b+b\cdot a)/2$ is a Jordan algebra denoted as $A^{(+)}$.
The natural question: what do all RB-operators on $M_2(F)^{(+)}$ look like, has not yet been considered. It is known that every Rota---Baxter operator of weight~$\lambda$ on an associative algebra~$A$ is an RB-operator of the same weight on $A^{(+)}$.
In 2024, V.N. Zhelyabin asked whether all RB-operators of weight zero on $M_2(F)^{(+)}$ come from $M_2(F)$. In 2024, he found a counterexample, it is an operator~$R_1$, which acts in terms of the matrix unities as follows:
$$
R_1(e_{11}) = e_{12}, \quad 
R_1(e_{12}) = 0, \quad 
R_1(e_{21}) = 0, \quad 
R_1(e_{22}) = - e_{12}.
$$
Hence, we may ask if there are other RB-operators of weight zero on $M_2(F)^{(+)}$ that are not RB-operators of weight zero on $M_2(F)$.
We may state the same question about the nonzero weight case.

The current work is devoted to the solution of all the questions mentioned above.

We find out that there are, up to all natural equivalences that preserve Rota---Baxter operators, exactly two RB-operators of weight zero on $M_2(F)^{(+)}$ that are not RB-operators of weight zero on $M_2(F)$: the one constructed by V.N. Zhelyabin and the following one:
$$
R_2(e_{11}) = 0, \quad 
R_2(e_{12}) = e_{11}, \quad 
R_2(e_{21}) = e_{11}, \quad 
R_2(e_{22}) = 0.
$$
Note that for derivations, the formal inverses of Rota---Baxter operators, one has a~coincidence of them on $M_n(F)^{(+)}$ and $M_n(F)$~\cite{DerJordan}.

Further, we check that both operators $R_1$ and $R_2$ satisfy the operator relation
$$
2R(x)R(y) 
 = R( R(x)y + xR(y) + yR(x) + R(y)x + \lambda xy + \lambda yx ) 
$$
for $\lambda = 0$.
We call an operator~$R$ satisfying this relation a~symmetrized RB-operator of weight~$\lambda$.

The reason why symmetrized RB-operators arise is simple.
Every symmetrized RB-operator of weight~$\lambda$ on an algebra~$A$ is
an RB-operator of weight~$\lambda$ on $A^{(+)}$.

We prove that every Rota---Baxter operator of weight~0 on $M_2(F)^{(+)}$
is an ordinary RB-operator or a symmetrized RB-operator of weight~0 on $M_2(F)$.

In the case of nonzero weight, the situation is slightly different.
The only RB-operators of nonzero weight on $M_2(F)^{(+)}$
that are not RB-operators of nonzero weight on $M_2(F)$
are splitting ones corresponding to decompositions $M_2(F)^{(+)} = H_2(F) \oplus B$ for some one-dimensional subalgebra $B$,
where $H_2(F)$ is, up to the action of automorphisms of $M_2(F)^{(+)}$, the only 
subalgebra of $M_2(F)^{(+)}$, which is not a subalgebra of $M_2(F)$. $H_2(F)$ is a Jordan subalgebra of symmetric matrices of order two.
We prove that every Rota---Baxter operator of weight~1 on $M_2(F)^{(+)}$
is an ordinary RB-operator of weight~1 on $M_2(F)$ or, up to the action of $\phi\colon R\to -R-\lambda\id$, is a symmetrized RB-operator of weight~1 on $M_2(F)$.
Since $H_2(F)$ is not associative, we may not avoid this action of~$\phi$.

We apply computer algebra~\texttt{Singular}~\cite{Singular} to decrease the routine with equations on parameters.
Throughout the work we suppose that the characteristic of the ground field~$F$ is not two.

\section{Preliminaries}

\subsection{Basic properties of Rota---Baxter operators}

We write down some basic properties of Rota---Baxter operators which can be found in~\cite{BGP,GuoMonograph}.

{\bf Proposition 1}.
Let $A$ be an algebra, and let $R$ be an RB-operator of weight~0 on~$A$.
Then the operator $\alpha R$ for any $\alpha\in F$ is again an RB-operator of weight~0 on~$A$.

{\bf Proposition 2}.
Given an algebra $A$, an RB-operator $R$ of weight $\lambda$ on $A$,
and $\psi\in\Aut(A)$, the operator $R^{(\psi)} = \psi^{-1}R\psi$
is again an RB-operator of weight $\lambda$ on $A$.

The same result is true when $\psi$ is an antiautomorphism of $A$, which means an automorphism of the vector space~$(A,+)$
such that $\psi(xy) = \psi(y)\psi(x)$ for all $x,y\in A$.
Transpose on the matrix algebra gives an example of an antiautomorphism of $M_n(F)$.

{\bf Proposition 3}.
Let $A$ be a unital algebra, and let $R$ be an RB-operator of weight~0 on~$A$.
Then $1\not \in \Imm(R)$.
Moreover, if $A$ is a~simple finite-dimensional algebra, $\dim A>1$, then $\dim \ker R\geq2$.

Given an algebra $A$ with a product $\cdot$, define the operations $\circ$ and $[\,,]$
on the vector space of $A$ by the rule
$$
a\circ b = \frac{a\cdot b+b\cdot a}{2},\quad 
[a,b] = a\cdot b-b\cdot a.
$$
We denote the space $A$ with $\circ$ as $A^{(+)}$
and the space $A$ with $[\,,]$ as $A^{(-)}$.

{\bf Proposition 4}.
Given an RB-operator $R$ of weight $\lambda$ on an algebra $A$,
$R$ is also an RB-operator of weight $\lambda$ on both $A^{(+)}$ and $A^{(-)}$.

Given an algebra $A$, let us define a map $\phi$ on the set of all RB-operators on $A$
as $\phi(P)=-P-\lambda(P)\id$, where $\lambda(P)$ denotes the weight of an RB-operator~$P$.
It is clear that $\phi^2$ coincides with the identity map.

{\bf Proposition 5}.
Let an algebra $A$ split as a vector space
into a direct sum of two subalgebras $A_1$ and $A_2$.
An operator $P$ defined as follows,
\begin{equation}\label{Split}
P(a_1 + a_2) = -\lambda a_2,\quad a_1\in A_1,\ a_2\in A_2,
\end{equation}
is an RB-operator of weight~$\lambda$ on~$A$.

An RB-operator from Proposition~5 is called a~splitting RB-operator with subalgebras $A_1,A_2$.
Note that the set of all splitting RB-operators on
an algebra $A$ is in bijection with all decompositions of~$A$
into a~direct sum of two subalgebras $A_1,A_2$.

{\bf Proposition 6}~\cite{BGP}.
Let $A$ be a unital algebra, and let $R$ be an RB-operator of weight~1 on~$A$.
If $R(1) \in F1$, then $R$ is splitting.

We call RB-operators $R = 0$ and $R = -\lambda\id$ trivial ones.

\subsection{Rota---Baxter operators on $M_2(F)$}

Let us list down all Rota---Baxter operators of weight zero on the matrix algebra of order two.
In~\cite{Aguiar} it was stated without proof over~$\mathbb{C}$.
In~\cite{Mat2}, this description was given with proof but without conjugation with automorphisms.
Next, in~\cite{BGP} this version appears over an algebraically closed field~$F$ of characteristic not two.
Finally, in~\cite{Panasenko}, the author avoids the condition that $F$ is algebraically closed.

{\bf Theorem 1}~\cite{Aguiar,BGP,Panasenko,Mat2}.
All nonzero RB-operators of weight zero on $M_2(F)$ over a~field~$F$ 
of characteristic not two up to conjugation with automorphisms of $M_2(F)$, transpose
and multiplication on a nonzero scalar are the following:

(A1) $R(e_{21}) = e_{12}$, $R(e_{11}) = R(e_{12}) = R(e_{22}) = 0$;

(A2) $R(e_{21}) = e_{11}$, $R(e_{11}) = R(e_{12}) = R(e_{22}) = 0$;

(A3) $R(e_{21}) = e_{11}$, $R(e_{22}) = e_{12}$, $R(e_{11}) = R(e_{12}) = 0$;

(A4) $R(e_{21}) =- e_{11}$, $R(e_{11}) = e_{12}$, $R(e_{12}) = R(e_{22}) = 0$. \\
\noindent Moreover, these 4 cases lie in different orbits
of the set of RB-operators of weight~0 on $M_2(F)$
under multiplication on a nonzero scalar or conjugation with $\Aut(M_2(F))$ or transpose.

Now, we provide the classification of RB-operators of nonzero weight on~$M_2(F)$.
All non-splitting RB-operators of nonzero weight were described in~\cite{BGP},
the final version over an algebraically closed field~$F $ of characteristic zero appeared in~\cite{GonGub}.

{\bf Theorem 2}~\cite{BGP,GonGub}.
Let $F$ be an algebraically closed field of characteristic zero.
Every nontrivial RB-operator on $M_2(F)$ of weight~1,
up to conjugation with an automorphism of $M_2(F)$ or transpose
and up to the action of $\phi$, equals one of the following cases:

(B1) $R(e_{11}) = e_{22}$, $R(e_{12}) = -e_{12}$, $R(e_{21}) = R(e_{22}) = 0$,

(B2) $R(e_{11}) = -e_{11}$, $R(e_{12}) = -e_{12}$, $R(e_{21}) = R(e_{22}) = 0$,

(B3) $R(e_{21}) = -e_{21}$, $R(e_{11}) = R(e_{12}) = R(e_{22}) = 0$,

(B4) $R(e_{21}) = e_{11}-e_{21}$, $R(e_{11}) = R(e_{12}) = R(e_{22}) = 0$,

(B5) $R(e_{12}) = e_{11} - e_{12}$, $R(e_{21}) = e_{22}-e_{21}$, $R(e_{11}) = R(e_{22}) = 0$,

(B6) $R(e_{11}) = -e_{11}$, $R(e_{12}) = -e_{12}$, $R(e_{21}) = 0$, $R(e_{22}) = e_{11}$.

\noindent Moreover, these 6 cases lie in different orbits
of the set of RB-operators of weight~1 on $M_2(F)$
under the action of $\phi$ and conjugation with $\Aut(M_2(F))$ or transpose.

The only non-splitting RB-operator of nonzero weight on~$M_2(F)$ is~(B1).

\subsection{Subalgebras of $M_2(F)^{(+)}$}

Let $\varphi\in\Aut(M_2(F)^{(+)})$, then by the Herstein theorem~\cite{HersteinTheorem}
one has that $\varphi$ is either an automorphism of $M_2(F)$
or an antiautomorphism of $M_2(F)$. 

Let $A$ be a nonzero proper subalgebra of $M_2(F)^{(+)}$.

{\sc Case 1}: $\dim A = 1$.
Hence, $A = L(a)$ for a nonzero~$a$.
Therefore $A$ is an associative subalgebra of $M_2(F)$.
Depending on whether $a^2 = 0$ or $a^2 = \lambda a$ for some $\lambda \neq 0$,
we get, up to conjugation, the following three variants: 
$A = Fe_{11}$,
$A = FE$,
$A = Fe_{12}$.

{\sc Case 2}: $\dim A = 2$.

It is easy to show that $A$ has to contain a non-nilpotent matrix. 

First, we consider the subcase when $A$ contains only degenerate matrices.
By the previous remark, it implies that $A$ contains an idempotent of rank~1.
Up to conjugation, we may assume that $e_{11}\in A$.
Consider $v =\begin{pmatrix}
0 & b \\
c & d \\
\end{pmatrix}$ such that $A = L(e_{11},v)$.
Since $e_{11}\circ v \in A$, we conclude that $d = 0$.
Further, $v^2 = bcE \in A$, thus either $b = 0$ or $c = 0$.
Up to conjugation with transpose, $A = L(e_{11},e_{12})$.

Second, we study the subcase when $A = L(E,v)$.
Since every matrix of order two satisfies the Cayley---Hamilton theorem,
we have $v^2 \in A$ for any non-scalar matrix~$v$.
Assuming that $F$ is quadratically closed, we conjugate~$v$ to its Jordan normal form. Hence, we have the following two subcases:
$A = L(e_{11},e_{22})$ and $A = L(E,e_{12})$.

{\sc Case 3}: $\dim A = 3$.
It is known that a maximal dimension of a subspace of $M_2(F)$ consisting only of degenerate matrices is two. 
Hence, $E \in A$.
Suppose that $E$ is the unique idempotent in~$A$.
Then we may find a nilpotent matrix in~$A$
and without loss of generality, consider the case $E,e_{12}\in A$.
Let $v$ be such that $A = L(E,e_{12},v)$.
We may assume that $v = \begin{pmatrix}
a & 0 \\
c & -a
\end{pmatrix}$. 
If $a = 0$, then $e_{21}\in A$ and $(e_{12}+e_{21})/\sqrt{2}$ is an idempotent matrix from~$A$, a contradiction. If $a \neq 0$, then we may assume that
$v = \begin{pmatrix}
1 & 0 \\
c & -1
\end{pmatrix}$.
If $c = 0$, then $e_{11},e_{22}\in A$, a contradiction.
Otherwise, define $u = (\sqrt{2/c})v$ and note that $u^2 = 1$.
In this case we again find two linear idempotent idempotents~$(u \pm E)/2$ in~$A$,
a~contradiction.

So up to conjugation, $A = L(e_{11},e_{22},v)$ for some matrix $v$ of the form
$be_{12} + ce_{21}$.
If $bc = 0$, then up to conjugation with transpose we have $A = L(e_{11},e_{22},e_{12})$.
Otherwise, we divide $v$ on $\sqrt{bc}$ and get the matrix $v' = \lambda e_{12} + (1/\lambda)e_{21}$.
Let us introduce the following automorphism $\psi_r$ of $M_2(F)$ for $r\in F\setminus \{ 0 \}$:
\begin{equation} \label{psi}
\begin{gathered}
\psi_r(e_{ii}) = e_{ii}, \ i=1,2,\quad
\psi_r(e_{12}) = re_{12}, \quad
\psi_r(e_{21}) = (1/r)e_{21}.
\end{gathered}
\end{equation}
Application of $\psi_{1/\lambda}$ gives us $A = L(e_{11},e_{22},e_{12}+e_{21})$.

{\bf Proposition 7}.
Up to the action of $\Aut(M_2(F)^{(+)})$, the list of all nonzero proper subalgebras of $M_2(F)^{(+)}$ is the following:
$Fe_{11},Fe_{12},FE,L(e_{11},e_{12}),L(e_{11},e_{22}),L(E,e_{12})$, \linebreak
$L(e_{11},e_{22},e_{12}), L(e_{11},e_{22},e_{12}+e_{21})$.

Note that $H_2(F) = L(e_{11},e_{22},e_{12}+e_{21})$ is, up to the action of $\Aut(M_2(F))^{(+)}$, the only subalgebra of $M_2(F)^{(+)}$ that is not a~subalgebra of $M_2(F)$.

\section{Case of weight zero}

Let $R$ be a nonzero Rota---Baxter operator of weight~0 on $M_2(F)^{(+)}$.
By Proposition~3, $\Imm R$ consists of only degenerate matrices, $\dim \Imm R\leq 2$,
we derive from the list of subalgebras of $M_2(F)^{(+)}$
that either $\Imm R = L(e_{11})$, or $\Imm R = L(e_{12})$,
or $\Imm R = L(e_{11},e_{12})$.
Anyway, we may assume that projections of $\Imm R$ on $e_{21}$ and $e_{22}$ are zero.
Denote
\begin{equation} \label{Interform}
\begin{gathered}
R(e_{11}) = a_{11}e_{11}+a_{12}e_{12}, \quad
R(e_{12}) = b_{11}e_{11}+b_{12}e_{12}, \\
R(e_{21}) = c_{11}e_{11}+c_{12}e_{12}, \quad
R(e_{22}) = d_{11}e_{11}+d_{12}e_{12}.
\end{gathered}
\end{equation}

Computations with \texttt{Singular} imply that $a_{11} = - b_{12} = -d_{11}$ and the following system of equations on the parameters holds:
\begin{gather*}
d_{11}^2 +b_{11}a_{12} = 0, \quad
b_{11}c_{12} = c_{11}d_{11}, \quad
d_{11}c_{12} + a_{12}(d_{12}+c_{11}+a_{12}) = 0, \\
a_{11}(d_{12}+a_{12}) = 0, \quad
b_{11}(d_{12}+a_{12}) = 0, \\
c_{12}(d_{12}+a_{12}) = 0, \quad
(d_{12}+a_{12})(d_{12}-a_{12}-c_{11}) = 0.
\end{gather*}

{\sc Case 1}: $a_{12} = 0$.
Then 
$$
e_{11}\in \ker R, \quad 
e_{12} \to b_{11} e_{11}, \quad
e_{21} \to c_{11} e_{11} + c_{12} e_{12}, \quad
e_{22} \to d_{12} e_{12},
$$
and the following equations hold:
$$
b_{11}c_{12} = 0, \quad
b_{11}d_{12} = 0, \quad
d_{12}c_{12} = 0, \quad
d_{12}(d_{12}-c_{11}) = 0.
$$
If $b_{11} \neq 0$, then $c_{12} = d_{12} = 0$, and $\dim (\Imm R) = 1$, 
we will consider this case below.

Let $b_{11} = 0$.
If $d_{12} = 0$, then we obtain again $\dim (\Imm R) = 1$, 
we will consider this case below.
If $d_{12} \neq 0$, then $c_{12} = 0$ and $c_{11} = d_{12}$.
Dividing on $d_{12}$, we obtain the operator
$$
e_{11},e_{12} \to 0, \quad 
e_{21} \to e_{11}, \quad 
e_{22} \to e_{12}.
$$

{\sc Case 2}: $a_{12} \neq 0$.
Dividing on~$a_{12}$, we may assume that $a_{12} = 1$.

If $d_{12} = -1$, then 
$$
R(e_{11}) = -d_{11}e_{11} + e_{12}, \quad
R(e_{12}) = d_{11}R(e_{11}), \quad
R(e_{21}) = c_{12}R(e_{11}), \quad
R(e_{22}) = -R(e_{11}),
$$
hence, $\dim (\Imm R) = 1$, see this case below.

If $d_{12}\neq-1$, we derive that $a_{11} = b_{11} = c_{12} = d_{12} = 0$,
$c_{11} = -1$, it is the operator
$$
e_{11} \to e_{12}, \quad 
e_{12},e_{22} \to 0, \quad
e_{21} \to -e_{11}.
$$

Finally, we study the case $\dim(\Imm R) = 1$.
If $\Imm R = L(e_{11})$, then \texttt{Singular} gives us the operator
$$
e_{11}\to0, \quad 
e_{12} \to b_{11} e_{11}, \quad
e_{21} \to c_{11} e_{11}, \quad 
e_{22} \to 0.
$$

If $\Imm R = L(e_{12})$, then computations imply
$$
e_{11} \to a_{12} e_{12}, \quad 
e_{12} \to 0, \quad
e_{21} \to c_{12} e_{12}, \quad 
e_{22} \to -a_{12} e_{12}.
$$

Below, we have the complete list of RB-operators of weight zero on~$M_2(F)^{(+)}$

1. $e_{11} \to e_{12}$, $e_{12} \to 0$,
$e_{21} \to - e_{11}$, $e_{22} \to 0$;

2. $e_{11},e_{12} \to 0$, $e_{21} \to e_{11}$, $e_{22} \to e_{12}$;

3. $e_{11}\to0$, $e_{12} \to \alpha e_{11}$,
$e_{21} \to \beta e_{11}$, $e_{22} \to 0$;

4. $e_{11} \to \alpha e_{12}$, $e_{12} \to 0$,
$e_{21} \to \beta e_{12}$, $e_{22} \to -\alpha e_{12}$.

Let us consider the third case $R(e_{11})=0$, $R(e_{12})=\alpha e_{11}$, $R(e_{21})=\beta e_{11}$, $R(e_{22})=0$, $(\alpha,\beta)\neq (0,0)$. 
Up to conjugation with transpose and up to multiplication by a scalar, we can assume that $\alpha = 1$. 
If $\beta\neq 0$, then the automorphism $\psi_{\sqrt{\beta}}$ is well-defined. 
If $R_1 = \frac{1}{\sqrt{\beta}}\psi_{\sqrt{\beta}}^{-1}R\psi_{\sqrt{\beta}}$, then 
$R_1(e_{11})=R_1(e_{22})=0$, $R_1(e_{12}) = e_{11}$, $R_1(e_{21})=e_{11}$. 

Let us consider the fourth case $R(e_{11})=\alpha e_{12}$, $R(e_{12})=0$, $R(e_{21})=\beta e_{12}$, $R(e_{22})=-\alpha e_{12}$. 
If $\alpha = 0$, then dividing on $\beta$ we have $R(e_{11})=R(e_{22})=R(e_{12})=0$ and $R(e_{21})=e_{12}$. 
If $\alpha\neq 0$, then we can assume that $\alpha = 1$. 
If $\beta\neq 0$ and $R_1 = \beta\psi_{\beta}^{-1}R\psi_{\beta}$, then $R_1(e_{11})=R_1(e_{21})=-R_1(e_{22})=e_{12}$, $R_1(e_{12})=0$. Let us define an automorphism $\varphi_t$ of $M_2(F)^{(+)}$ as follows:
\begin{equation} \label{varphi_t}
\begin{gathered}
\varphi_t(e_{11})=e_{11}+t e_{12}, \quad
\varphi_t(e_{12})=e_{12}, \\ 
\varphi_t(e_{21}) = e_{21} - t e_{11} + t e_{22} - t^2 e_{12}, \quad 
\varphi_t(e_{22})=e_{22}-t e_{12}.
\end{gathered}
\end{equation}
For $R_2=\varphi_{1/2}^{-1}R_1\varphi_{1/2}$, we have 
$R_2(e_{11})=-R(e_{22})=e_{12}$, $R(e_{12})=R(e_{21})=0$, so we have a~case $\beta = 0$.

We have described all RB-operators of weight zero on $M_2(F)^{(+)}$.

{\bf Theorem 3}.
Every nonzero RB-operator of weight zero on $M_2(F)^{(+)}$ over a quadratically closed~field~$F$ 
of characteristic not two up to conjugation with automorphisms of $M_2(F)$, transpose
and multiplication on a nonzero scalar either coincides with one of the operators~(A1)--(A4) or is one of the following:

(R1) $R(e_{11}) = e_{12}$, $R(e_{12}) = 0$, $R(e_{21}) = 0$, $R(e_{22}) = - e_{12}$;

(R2) $R(e_{11}) = 0$, $R(e_{12}) = e_{11}$, $R(e_{21}) = e_{11}$, $R(e_{22}) = 0$.

\noindent Moreover, these 6 cases lie in different orbits
of the set of RB-operators of weight~0 on $M_2(F)^{(+)}$
under multiplication on a nonzero scalar or conjugation with $\Aut(M_2(F)^{(+)})$.

{\sc Proof}.
It remains to prove the second part of Theorem.
Since all actions listed in the formulation preserve the set of RB-operators on $M_2(F)$,
the operators~(A1)--(A4), on the one hand, and the operators (R1) and (R2), on the other hand,
lie in different orbits. 
Indeed,
$$
e_{11} = R_1(e_{12})R_1(e_{21}) \neq 0, \quad
0 = R_2(e_{21})R_2(e_{11}) \neq R_2(e_{22}) = -e_{12}.
$$

By Theorem~1, all operators~(A1)--(A4) lie in different orbits.
We distinguish the operators~(R1) and~(R2), since their images are not (anti)isomorphic subalgebras.
\hfill $\square$

\section{Case of nonzero weight}

{\bf Theorem 4}.
Every non-trivial RB-operator of weight~1 on $M_2(F)^{(+)}$ over a~quadratically closed field~$F$ 
of characteristic not two up to conjugation with automorphisms of $M_2(F)$, transpose
and multiplication on a nonzero scalar either coincides with one of the operators~(B1)--(B6) or 
is one of the following:

(S1) $R(e_{11}) = R(e_{22}) = R(e_{12}+e_{21}) = 0$, $R(e_{12}) = -e_{12}$;

(S2) $R(e_{11}) = R(e_{22}) = R(e_{12}+e_{21}) = 0$, $R(e_{12}) = -e_{11}-e_{12}$.

\noindent Moreover, these 8 cases lie in different orbits
of the set of RB-operators of weight~1 on $M_2(F)^{(+)}$
under multiplication on a nonzero scalar or conjugation with $\Aut(M_2(F)^{(+)})$.

{\sc Proof}.
By~\cite{Spectrum}, we have that $R(1)$, up to the action of $\phi$, is either zero or conjugate either to $e_{22}$ or to $-e_{11}$.

{\sc Case 1}: $R(1) = e_{22}$.
Computations with \texttt{Singular} provide us with a pair of RB-operators that are conjugate under transpose. Thus, we write down only one of them:
$$
e_{11} \to e_{22}, \quad 
e_{12}, e_{22} \to 0, \quad 
e_{21} \to -e_{21}.
$$

{\sc Case 2}: $R(1) = -e_{11}$.
Again, \texttt{Singular} gives us a pair of RB-operators that are conjugate under transpose. Below, we have one of them:
$$
e_{11} \to -e_{11}, \quad e_{12} \to -e_{12}, \quad
e_{21}, e_{22} \to 0.
$$

{\sc Case 3}: $R(1) = 0$, hence, $R$ is splitting by Proposition~6.
Denote $M_2(F)^{(+)} = A\oplus B$.
If $\dim A = \dim B = 2$, then by Proposition~7 we have a decomposition onto two subalgebras of~$M_2(F)$. Hence, such an operator is an RB-operator on~$M_2(F)$.
Now, we consider the case $\dim A = 3$ and $\dim B = 1$.
If $A$ is isomorphic to $L(e_{11},e_{12},e_{22})$,
then both $A$ and $B$ are associative, and again $R$ is an RB-operator on~$M_2(F)$, all such decompositions are listed in Theorem~2.

Finally, we consider $A = L(e_{11},e_{22},e_{12}+e_{21})$, then 
$B = L(v)$
for $v = \begin{pmatrix}
x & y \\
z & t 
\end{pmatrix}$, where $y\neq z$.
Suppose that $v$ is nilpotent, hence $t = -x$ and $x^2 = yz = 0$.
If $x = 0$, then $v \in L(e_{12})$ or $v\in L(e_{21})$.
Up to conjugation with transpose, we obtain the decomposition $L(e_{11},e_{22},e_{12}+e_{21}) \oplus L(e_{12})$.
If $x\neq0$, then $y,z\neq0$. Dividing on~$z$, we may suppose that $v$ has the form
$v = \begin{pmatrix}
-t & -t^2 \\
1 & t
\end{pmatrix}$. Since $v\not \in A$, we have $t^2+1 \neq 0$.
Note that~$A$ is invariant under conjugation with a~matrix
$P_a = \begin{pmatrix}
a & \sqrt{1-a^2} \\
- \sqrt{1-a^2} & a
\end{pmatrix}$, where $a\in F$.
We apply conjugation with $P_{t/\sqrt{t^2+1}}$ and get $v \in L(e_{12})$.

Suppose that $v$ is an idempotent, that is $x = 1 - t$, $(1-t)t - yz =0$.
If $y = 0$ or $z = 0$, then up to conjugation with an automorphism of $M_2(F)^{(+)}$,
we obtain the decomposition $L(e_{11},e_{22},e_{12}+e_{21}) \oplus L(e_{11}+e_{12})$.
Otherwise, we express $v = \begin{pmatrix}
1-t & y \\
\frac{(1-t)t}{y} & t 
\end{pmatrix}$, moreover, $Q = \frac{y^2 + (1-t)t}{y}\neq0$, since $v\not \in A$.
We want to conjugate $v$ by $P_a$ to make the image of $v$ with a~zero coordinate on $e_{11}$ or on $e_{22}$. Thus, we have to find~$a$ such that
$$
a^2(1-2t) - a\sqrt{1-a^2}Q = \alpha - t, 
$$
for at least one of the values $\alpha\in \{0,1\}$.
If $t = 1/2$, then we find~$a$ from the equation $a^2(1-a^2) = (\alpha-t)^2/Q^2$.
Otherwise, squaring the expression implies
$$
a^4((2t-1)^2 + Q^2) + a^2( 2(\alpha-t)(2t-1) - Q^2 ) + (\alpha-t)^2 = 0.
$$
We may resolve it for~$a$ if at least one of the coefficients at $a^4$ and $a^2$ is nonzero.
Suppose to the contrary, that both these coefficients are zero.
Hence, $Q^2 = -(2t-1)^2 = 2(\alpha-t)(2t-1)$.
Since $t\neq1/2$, we get $1-2t = 2\alpha-2t$, which does not hold for $\alpha\in \{0,1\}$.
Therefore, conjugation with~$P_a$ for appropriate~$a$ sends $v$ to a vector from
$L(e_{11}+e_{12})$.

The proof of the second part of Theorem is analogous to the proof of~Theorem~3.
\hfill $\square$

{\bf Remark 1}.
Note that unlike the Lie case, one may decompose a simple Jordan algebra as a direct vector space sum of two simple ones: $M_2(F)^{(+)} \cong H_2(F) \oplus F$, the operator~(S1) corresponds to this decomposition.
For more about simple decompositions of simple Jordan algebras, see in~\cite{Tvalavadze}.

Let us introduce a new kind of operator.

{\bf Definition 2}.
Let $A$ be an algebra, a linear operator $R$ on $A$ is called symmetrized Rota---Baxter operator (symmetrized RB-operator, for short) of weight~$\lambda$, if the following relation
\begin{equation} \label{half-RB}
2R(x)R(y) 
 = R( R(x)y + xR(y) + yR(x) + R(y)x + \lambda xy + \lambda yx ) 
\end{equation}
holds for all $x,y \in A$.

{\bf Proposition 8}.
Given a symmetrized RB-operator $R$ of weight $\lambda$ on an algebra $A$,
$R$ is an RB-operator of weight~$\lambda$ on $A^{(+)}$.

{\sc Proof}.
We have
\begin{multline*}
2 R(x)\circ R(y)    
 = R(x)R(y) + R(y)R(x) \\
 = R( R(x)y + yR(x) + R(y)x + xR(y) + \lambda xy + \lambda yx) \\
 \hfill = 2R( R(x)\circ y + x\circ R(y) + \lambda x\circ y).
 \hfill \square
\end{multline*}

{\bf Corollary 1}.
Let $A$ be an associative algebra, $R$ be a symmetrized RB-operator of weight~$\lambda$.
Then $(A,*)$, where 
\begin{equation} \label{*-product}
x*y = (R(x)y + xR(y) + yR(x) + R(y)x + \lambda xy + \lambda yx)/2, 
\end{equation}
is a Jordan algebra.
Indeed, $R$ is an RB-operator of weight~$\lambda$ on the special Jordan algebra $(A^{(+)},\circ)$.
It is known~\cite{preJordan} that the operation 
$x\star y = R(x)\circ y + x\circ R(y) + \lambda x\circ y$ defines a Jordan product on the space~$A$. 
Here we have $x\star y = x*y$.

Define polynomials
$$
H_{r,s}(x) = (x-r)(x-r+1)\ldots x(x+1)\ldots (x+s-1)
$$
for natural $r,s\geq0$.

{\bf Corollary 2}.
Let $A$ be a unital finite-dimensional algebra and let $R$ be a symmetrized RB-operator of weight~$\lambda$ on~$A$.

a) Then there exist $k,l\geq0$ such that $R^k(R+\lambda\id)^l = 0$. 
Hence, $\Spec(R)\subset\{0,-\lambda\}$;

b) If $\lambda = 0$, then $R(1)$ is nilpotent;

c) If $\lambda = -1$, then there exist $r,s$ such that the minimal polynomial of $R(1)$ equals $H_{r,s}(x)$.

{\sc Proof}.
By~Proposition~8, $R$ is an RB-operator of weight~$\lambda$ on~$A^{(+)}$,
hence, we may apply the results from~\cite{Spectrum}, since $A^{(+)}$ is again unital.
\hfill $\square$

Let us note some basic properties of such operators.
Surely, $\Imm R$ is a subalgebra of~$A$.
The analogues of Propositions~1--3 hold for symmetrized RB-operators too.
For example, to prove the analogue of Proposition~3 it is enough to 
consider a symmetrized RB-operator $R$ of weight~0 defined on an algebra $A$
as an RB-operator of weight~0 on~$A^{(+)}$.

Since the right part of~\eqref{half-RB} is symmetric with respect to $x$ and $y$, so we have $R(x)R(y) = R(y)R(x)$ for any $x,y\in A$. 
It means that $\Imm R$ is a commutative subalgebra of~$A$.
It implies that the set of all symmetrized Rota---Baxter operators on a given algebra~$A$ is closed under conjugation with an antiautomorphism of~$A$. 

The lack of symmetry is not the only reason why, as far as we know, this kind of operator has not yet been studied.
If one tries to interpret~\eqref{half-RB} as a homomorphism condition $R\colon (A,*)\to (A,\cdot)$,
where $*$ is defined by~\eqref{*-product}, then the new operation~$*$ may not satisfy the properties that are for the initial product~$\cdot$. 
For example, if $\cdot$ is associative, $*$ may not be.
It is a reason why symmetrized RB-operators do not appear in the list of, in some sense, good operators in~\cite{RotaProblem}.

If an algebra $A$ is associative and commutative, then the notions of ordinary and symmetrized Rota---Baxter operators coincide.
If $A$ is anticommutative, then a symmetrized RB-operator~$R$ is a linear operator satisfying the identity $R(x)R(y) = 0$, in~\cite{Gub2016} such operators were called abelian.

Let us return to RB-operators on $M_2(F)^{(+)}$ and formulate the deep connection between them and symmetrized RB-operators on $M_2(F)$.

{\bf Corollary 3}.
Let $F$ be a quadratically closed field of characteristic not two. Then

a) every RB-operator of weight~0 on $M_2(F)^{(+)}$
is either an RB-operator of weight~0 on $M_2(F)$
or a symmetrized RB-operator of weight~0 on $M_2(F)$;

b) every RB-operator of weight~1 on $M_2(F)^{(+)}$
is either an RB-operator of weight~1 on $M_2(F)$
or, up to the action of $\phi$, a symmetrized RB-operator of weight~1 on $M_2(F)$.

Note that a linear operator may be simultaneously an ordinary and a symmetrized RB-operator, e.\,g. the operators (A1) and (B3) are so. 

{\bf Remark 2}.
It looks like one can not adopt Corollary~3 for the case of $M_2(F)^{(-)}$.
We may consider the known description of Rota---Baxter operators of weight~1 on $\mathrm{gl}_2(F)$~\cite{Goncharov2022} and check if all operators from the list satisfies
the relation of the form
$R(x)R(y) = R(\Sigma)$,
where $\Sigma$ is a linear combination of $R(x)y,yR(x),R(y)x,xR(y),xy,yx$.
Maybe, here we should involve operators of the form 
$R(x)R(y) = R(\Sigma_1) + \Sigma_2$.

\section*{Acknowledgements}

The authors are grateful to V.N. Zhelyabin and P.S. Kolesnikov for the helpful discussions.
The study was supported by a grant from the Russian Science Foundation №\,23-71-10005, https://rscf.ru/project/23-71-10005/

\noindent Vsevolod Gubarev \\
Alexander Panasenko \\
Novosibirsk State University \\
Pirogova str. 1, 630090 Novosibirsk, Russia \\
Sobolev Institute of Mathematics \\
Acad. Koptyug ave. 4, 630090 Novosibirsk, Russia \\
e-mail: wsewolod89@gmail.com

\end{document}